\documentclass[a4paper,11pt]{article}
\usepackage[latin1]{inputenc}
\usepackage[all] {xy}
\usepackage[dvips]{graphicx}
\usepackage{amssymb,amsthm,amsmath,psfrag}
\usepackage[dvips]{geometry}

\newcommand{\N}{\mathbb{N}}
\newcommand{\R}{\mathbb{R}}

\usepackage[latin1]{inputenc}
\usepackage[all] {xy}
\usepackage[dvips]{graphicx}
\usepackage{amssymb,amsthm,amsmath,psfrag}
\usepackage[dvips]{geometry}

\def\rc{\ar[r]^<(.5)>(.5)}

\def\gt{G\textnormal{-i.u.t.a.}}

\newtheorem{Theorem}{Theorem}
\newtheorem{Definition}{Definition}
\newtheorem{Remark}{Remark}
\newtheorem{Example}{Example}
\newtheorem{Lemma}{Lemma}
\newtheorem{Proposition}{Proposition}
\newtheorem{Corollary}{Corollary}

\title{Compact Lie algebras,  transversely Lie  foliations and fibrations}
\author{M. Tavares}

\date{}

\begin{document}

\maketitle

\tableofcontents

\section{Introduction}

An important theorem due to Tischler states that  a smooth closed
manifold fibers over the q-dimensional torus if and only if admits a
foliation given by q linearly independent closed forms of degree
one. Indeed, given $q$ linearly independent closed forms of degree
one, it is possible to find new closed forms arbitrarily close to
the originals and having {\it rational periods}. The foliation
defined by these new forms has all leaves compact and these leaves
are the fibers of a locally trivial fibration over the torus. All
this can be seen as embedded in the theory of Lie-foliations
\cite{Fedida}, \cite{Go}. Let $\mathfrak{g}$ be a Lie algebra of
dimension $q$ and let $M$ be a smooth closed manifold. We say that a
foliation $\mathcal{F}$ of codimension $q$ on $M$ is a {\it
$\mathfrak{g}$-Lie foliation} if, given structure constants
$\{c_{ij}^k\}$ of $\mathfrak{g}$, there exists a system $
\{\omega_1, \dots,\omega_q \}$ of linearly independent smooth forms
of degree one on $ M $ which satisfies $ d \omega_k = - \sum_{i <j}
c_{ij}^k \omega_i \wedge \omega_j$ and whose kernel is tangent to
$\mathcal{F}$. Let $G$ be a connected Lie group having
$\mathfrak{g}$ as Lie algebra. A foliation $\mathcal{F}$ on $M$ is a
Lie $\mathfrak{g}$-foliation if and only if $\mathcal{F}$ admits a
{\it Lie transverse structure modeled by $G$} (see \cite{Go}), that
is, $\mathcal{F}$ is defined by local submersions taking values on
$G$ with transition maps given by left translations on $G$. The
existence of Lie foliation does not guarantee the existence of a
global submersion of $M$ in $G$. The difficulty in obtaining such a
submersion from a Lie foliation $\mathcal{F}$ can be measured by the
{\it global holonomy group of $\mathcal{F}$ with respect to $G$},
which is a representation of the fundamental group of $M$ in $G$. In
his thesis, Fedida proves that, in fact, the structure of
$\mathcal{F}$ is strongly related to its global holonomy
(\cite{Fedida}, \cite{Go}). In this work, we obtain a generalization
of the Tischler fibration Theorem for $\mathfrak{g}$-Lie foliations,
when $\mathfrak{g}$ is a {\it compact} Lie algebra, in the sense
that $\mathfrak{g}$ is the Lie algebra of some compact Lie group:

\begin{Theorem}
\label{Theorem:A} Let $\mathfrak{g}$  be a compact Lie algebra,  $M$
a smooth closed manifold and  $\mathcal{F}$ a Lie
$\mathfrak{g}$-foliation on $M$. Let $G$ be a connected compact Lie
group having $\mathfrak{g}$ as Lie algebra. If the closure of the
global holonomy of $\mathcal F$ with respect to $G$ has Abelian Lie
algebra, then $M$ possesses a finite covering which fibers over $G$.
Indeed,  there is a finite covering of $M$ where the lifting of
$\mathcal{F}$ can be $C^0$-approximated by Lie
$\mathfrak{g}$-foliations  having all leaves compact.
\end{Theorem}

The proof of Theorem~\ref{Theorem:A} gives a
description of the structure of the foliation, and this allows us to
well understand the finite covering mentioned.
Combining our result with a result due to Caron and Carrière (\cite{carcar}),
we have

\begin{Corollary}\label{Corollary:i}
 Let $\mathfrak{g}$  be a semisimple
compact Lie algebra and let $\mathcal{F}$ be a (not necessarily
minimal) one-dimensional Lie $\mathfrak{g}$-foliation on a  smooth
closed manifold $M$. We have
\begin{enumerate} \item there is a finite covering of $M$ where
the lifting of $\mathcal{F}$
can be $C^0$-approximated by Lie $\mathfrak{g}$-foliations having all leaves compact;
                  \item If $G$ is the simply connected Lie group having $\mathfrak{g}$ as Lie algebra, then $M$ possesses a finite covering
                      which
                  fibers over $G$.
\end{enumerate}

\end{Corollary}

An important remark in the theory is that the closure  of the global
holonomy of a Lie foliation must be an Abelian group provided that
the fundamental group of $M$ is Abelian. In this work, we extend
this observation by proving that if the fundamental group of $M$ is
amenable and $G$ is compact, then the closure of the global holonomy
of $\mathcal{F}$ has Abelian Lie algebra.
 Thus, we can apply Theorem A to obtain the following:

\begin{Theorem}
\label{Theorem:B} Let $\mathfrak{g}$ be a compact Lie algebra and
$M$ a smooth closed manifold with amenable fundamental group. If $M$
admits a Lie $\mathfrak{g}$-foliation, then there is a finite
covering of $M$ where the lifted foliation can be $C^0$-approximated
by Lie $\mathfrak{g}$-foliations having all leaves compact. If  $G$
is a compact Lie group having $\mathfrak{g}$ as Lie algebra, then
$M$ possesses   a finite covering which
                  fibers over $G$.
\end{Theorem}

We also studied the particular case where the fundamental group of the ambient manifold has subexponential growth and
$\mathfrak{g}$ is not necessarily compact:

\begin{Theorem}
\label{Theorem:C}
 Let $M$ be a smooth closed manifold with
fundamental group of subexponential growth and suppose that $M$
admits a Lie $\mathfrak{g}$-foliation. We have:
\begin{enumerate} \item If \,  $G$ is a connected Lie group having $\mathfrak{g}$ as Lie algebra, then $G$ has polynomial volume growth;
                  \item If $\mathfrak{s} < \mathfrak{g}$ is a Levi factor of $\mathfrak{g}$ and $S$ is a connected Lie group having
                      $\mathfrak{s}$ as
                  Lie algebra, then $S$ is compact and $M$ possesses a finite covering which fibers over $S$.
\end{enumerate}
\end{Theorem}

In the semisimple case, the result above has the following statement:

\begin{Corollary} \label{Corollary:ii}
Let $M$ be a smooth closed manifold which admits a Lie
$\mathfrak{g}$-foliation where  $\mathfrak{g}$  is a semisimple Lie
algebra. If  $\mathfrak{g}$ is non-compact  then the fundamental group of $M$ has
exponential growth.  If the fundamental group of $M$ has
subexponential growth then any  connected Lie group $G$ having
$\mathfrak{g}$ as Lie algebra  is compact and $M$ possesses a finite
covering which fibers over $G$.
\end{Corollary}

The concept of Lie foliations admits a natural generalization: Let $G$ be a connected Lie group of dimension $n$
 and let $K$ be a connected closed Lie subgroup of $G$ of dimension $n-q$. Let $\mathfrak{g}$ and $\mathfrak{k}$ be
  the Lie algebras of $G$ and $K$, respectively. Let $\beta = \{X_1, \dots, X_q, X_{q +1},\dots, X_n\}$ be a basis of $\mathfrak{g}$
such that $X_{q +1}, \dots, X_n $ is a basis of $\mathfrak{k}$ and
denote by $\{c_{ij}^k\} $ the structure constants of $\beta $.
Suppose that a smooth manifold $M$ admits $n$ smooth forms of degree
one $\omega_1, \dots, \omega_q, \omega_{q +1} , \dots,  \omega_n$
satisfying $d \omega_k = - \sum_{i <j} c_{ij}^k \omega_i \wedge
\omega_j $. If the forms $\omega_1, \dots, \omega_q $ are linearly
independent at each point of $ M $, then they define a foliation $
\mathcal{F}$ of codimension $q$ on $ M $ and, for this case, we can
adapt Theorem A. Such a foliation admits a structure of {\it
transversely homogeneous foliation modeled by $ G/K$} (see
\cite{blu}, \cite{Go}), that is, $ \mathcal F $ is defined by local
submersions taking values on $ G/K $ with  transition maps given by
left translations on $G$. Similar structures are constructed, as a
global holonomy group (\cite{Go}, \cite{blu}). We should note,
however, that not every transversely homogeneous foliation is
obtained from global forms of degree one. Our preceding results
together with  results of Blumenthal (\cite{blu}, \cite{blu1}), give
the following fibration result for transversely homogeneous
foliations:

\begin{Theorem}
\label{Theorem:D}
 Let $G$ be a connected compact Lie group,  $K$
a closed and connected subgroup of $G$ and let $M$ be a smooth
closed manifold with amenable fundamental group. Given
a $G/K$-transversely homogeneous foliation $\mathcal{F}$ on $M$,  suppose one of
the following conditions is satisfied

\begin{enumerate} \item $H^1(M,K)$ is trivial;
                  \item the normal bundle of $\mathcal{F}$
                  is trivial and $G/K$ is a simply-connected
                  Riemannian homogeneous space of constant curvature.
\end{enumerate}
Then there is a finite covering of $M$ where the lifting of
$\mathcal{F}$ can be $C^0$-approximated by $G/K$-transversely
homogeneous foliations. Moreover, $M$ possesses a finite covering
which fibers over $G/K$.
\end{Theorem}

We stress that conditions (1) and (2) in Theorem~\ref{Theorem:D}
above are used in order to assure that the homogeneous transverse
structure is given by globally defined one-forms and can be replaced
by this hypothesis, which may also verified in some other
situations.

\begin{Corollary} Let $\mathcal{F}$ be an one-dimensional $SO(2q+1)/SO(2q)$-transversely homogeneous foliation with trivial
normal bundle on a smooth closed manifold $M$. Then there is a finite covering of $M$ where the lifting of
$\mathcal{F}$ can be $C^0$-approximated by $SO(2q+1)/SO(2q)$-foliations having all leaves compact.
If $M$ is orientable and $q=1$, then the finite covering mentioned is either the identity or a double covering of $M$.
\end{Corollary}

\section{Lie foliations}

Let us recall the formal definition and main properties we need concerning Lie foliations (see \cite{Go}, \cite{Fedida} for more):

\begin{Definition}{\rm Let $M$ be a smooth manifold  and $\mathfrak{g}$
a Lie algebra of dimension $q$. We say that a codimension $q$ foliation
 $\mathcal{F}$ on $M$ is a {\it Lie $\mathfrak{g}$-foliation} if
$\mathcal{F}$ is tangent to the kernel of a smooth form of degree one $\Omega$ on
$M$ with values in $\mathfrak{g}$ such that
\begin{enumerate} \item $\Omega(x): T_x M \to \mathfrak{g}$ is surjective,
for each $x \in M$;
                  \item $d \Omega + \frac{1}{2} [\Omega,\Omega] =0$.

\end{enumerate}}
\end{Definition}

Let $X_1, \dots, X_q$ be a basis of a Lie algebra $\mathfrak{g}$ with structure constants $\{c_{ij}^k\}$.
If there are linearly independent forms of degree one $\omega_1, \dots,
\omega_q$ on a smooth manifold $M$ which satisfy $d  \omega_k= - \sum_{i <j} c_{ij}^k \omega_i
\wedge \omega_j$, then $\Omega = \sum_{l=1}^q \omega_l X_l$
is a form of degree one on $M$ with values in $\mathfrak{g}$ satisfying the second item of Definition 1.
Conversely, let $\Omega$ be a smooth form of degree one on $M$ with values in $\mathfrak{g}$
satisfying (2). If $X_1, \dots,
X_q$ is a basis of $\mathfrak{g}$ and $\Omega = \sum_{l=1}^q \omega_l
X_l$, then $d  \omega_k= - \sum_{i <j} c_{ij}^k \omega_i \wedge \omega_j$.
 The form $\Omega$ is surjective in each point of $M$ if and only if
 the forms $\omega_1, \dots, \omega_q$ are linearly independent in each point of $M$.
 Moreover, the kernel of $\Omega$ coincides with the kernel of the system $\{ \omega_1, \dots, \omega_q\}$.
 Therefore, a codimension $q$ foliation $\mathcal{F}$ on $M$ is a Lie $\mathfrak{g}$-foliation
if and only if $\mathcal{F}$ is tangent to the kernel of a system
$\{\omega_1, \dots, \omega_q\}$ of smooth linearly independent forms of degree one on $M$ which satisfies
$$d \omega_k = - \sum_{i<j} c_{ij}^k \omega_i \wedge \omega_j.$$

\begin{Example}{\rm A Lie $\R^r$-foliation on a smooth manifold $M$ is a
foliation defined by the kernel of $r$ smooth linearly independent closed forms
of degree one on $M$.}
\end{Example}

\begin{Example}{\rm  Let $\pi: E \to M$ be a $G$-principal bundle.
If a connection form $\Omega$ one $E$ is flat, that is, has zero
curvature, then $\Omega$ satisfies $d \Omega= - \frac{1}{2} [\Omega,
\Omega]$ (see \cite{kob}, pg 78, 92-94). Therefore, $\Omega$ defines
a Lie $\mathfrak{g}$-foliation $\mathcal{F}$ on $E$ which is
transverse to the fibers.}
\end{Example}

We also recall the concept of Lie transverse structure for foliations:

\begin{Definition}\label{liefol}\rm{Let $\mathcal{F}$ be a codimension $q$ foliation on a differentiable manifold $M$ and
let $G$ be a $q$-dimensional Lie group. A Lie transverse structure of model $G$ for $\mathcal{F}$ is given by a family of submersions $\{f_i:U_i \to G\}_{i \in \Lambda}$
and a family of locally constant maps $g_{ij}: U_i \cap U_j \to \{\textnormal{left translations on }G\}$ such that
\begin{enumerate}\item $\{U_i\}_{i \in \Lambda}$ is an open cover of $M$;
                 \item the plaques of $\mathcal{F}$ in $U_i$  are given by $f_i = {\rm constant}$;
                 \item $f_i(x) = g_{ij}(x) \circ f_j(x)$, $\forall x \in U_i \cap U_j.$
\end{enumerate}}
\end{Definition}

The result bellow allow us to understand the relation between Lie $\mathfrak{g}$-foliations and
a foliation with a Lie transverse structure:

{\bf \flushleft  Darboux-Fedida's Theorem.}{ \it Let $M$ be a smooth manifold and let $\mathcal{F}=
\mathcal{F}_{\Omega}$ be a Lie $\mathfrak{g}$-foliation on $M$.
Let $G$ be a connected Lie group having $\mathfrak{g}$ as Lie algebra and denote by $\alpha$ its Maurer-Cartan
form. Denote by $\phi$ the natural action of $G$ on $M \times
G$ by left translations and let $p_1: M \times G \to M$ and $p_2:
M \times G \to G$ be the canonical projections. We have that:
\begin{enumerate} \item The form $\Theta$ on $M \times G$ with values in
$\mathfrak{g}$ defined by
                  $$\Theta_{(x,g)}(u,v) = \Omega(x)u - \alpha(g)v$$
                  induces on $M \times G$ a codimension $q$ foliation $\mathcal{G}$;
                  \item  $\mathcal{G}$ is invariant by $\phi$, that is,
                  for each $g \in G$, the diffeomorphism $\phi^g$
                  take leaves of $\mathcal{G}$ over leaves of $\mathcal{G}$.
                  \item The leaves of $\mathcal{G}$ are transversal
                  to the orbits of $\phi$.
                  \item  If $M'$ is a leaf of $\mathcal{G}$, the map
                  $p_1|_{M'} : M' \to M$ is a Galoisian covering.
                  The automorphism group of this covering is associated
                  with a homomorphism $h: \pi_1(M) \to G$ which satisfies
                  $$\alpha(p) = \phi^{h (\alpha) } (p),
                  \ \forall \alpha \in \pi_1(M) \textnormal{ and } p \in M'$$
                  and whose image is $h(\pi_1(M)) =
                  \{g \in G ; \ \phi^g(M') = M' \}$;
                  \item The map $p_2|_{M'} : M' \to G$ is a submersion;
                  \item $p_1^*(\Omega)|_{M'} \equiv p_2^*(\alpha)|_{M'}$.
\end{enumerate}}

Consider the notations used in the Darboux-Fedida's Theorem. The
group $\Gamma = h(\pi_1(M)) < G$ is called {\it global holonomy
group of $\mathcal{F}$ with respect to $G$}, or, for simplicity,
{\it holonomy global of $\mathcal{F}$}. We observe that, fixed the
form $\Omega$ and the Lie group $G$, the global holonomy group of
$\mathcal{F}$ with respect to $G$ is well defined up to a
conjugation by an element of $G$, because it also depends on the
choice of the leaf $M'$ of foliation $\mathcal{G}$. Throughout the
work, if $\mathcal{F} =\mathcal{F}_{\Omega}$ is a Lie $\mathfrak
{g}$-foliation on $M$ and $G$ is a connected Lie group which has
$\mathfrak{g}$ as Lie algebra, we denote by $\Gamma(\Omega,G)$ the
global holonomy group of $\mathcal{F}$ constructed according to
Darboux-Fedida' Theorem, after having chosen a leaf $M'$ of the
foliation $\mathcal{G}$. Still using the notations of the
Darboux-Fedida Theorem, consider the Galoisian covering $P_1 =
p_1|_{M'}: M' \to M $, the submersion $P_2 = p_2|_{M'}: M' \to G$
and the homomorphism o $h: \pi_1(M) \to G$ associated to $P_1$. We
have that
\begin{itemize} \item[i.] $P_2$ is equivariant by $h$, that is, $P_2(\alpha(x)) = h(\alpha) * P_2(x)$ ;
                \item[ii.]$P_1^*(\mathcal{F})$ is the foliation defined by $P_2$;
\end{itemize}
The data $(M',P_2,G,h)$ is a {\it $G$-development for $\mathcal{F}$} (see \cite{Go}).
If we choose another leaf $M''$ of $\mathcal{G}$ different of $M'$, we obtain
another $G$-development $(M'',\tilde{P}_2,\tilde{h},G)$ for $\mathcal{F}$.
The two developments are related as follows: $\tilde{h} = g * h *
g^{-1}$ e $\tilde{P}_2 = g * P_2$, for an element $g \in G$ such that $\phi^g(M') = \tilde{M}$.
Using a development, it is easy to obtain a Lie transverse structure of model $G$ for the foliation.
Conversely, let $\mathcal{F}$ be a foliation that admits a Lie transverse structure of model $G$.
If $\alpha_1, \cdots, \alpha_q$ are
linearly independent forms of degree one in $G$ invariant by left translations with structure constants
$\{c_{ij}^k\}$, then
we can use the local submersions that define $\mathcal{F}$ to obtain smooth linearly independent global
forms of degree one $\omega_1, \dots, \omega_q$ on $M$
which satisfy $d \omega_k = \sum_{i
<j} c_{ij}^k \omega_i \wedge \omega_j$.
Therefore, a foliation
$\mathcal{F}$ is a Lie $\mathfrak{g}$-foliation if and only if $\mathcal{F}$ admits
a $G$-development. Moreover, $\mathcal{F}$ is a Lie
$\mathfrak{g}$-foliation, if and only if $\mathcal{F}$ admits
a Lie transverse structure of model $G$.

We say that a Lie algebra $\mathfrak{g}$ is {\it compact}
if it is the Lie algebra of some compact connected Lie group.
The example bellow describes a known method of obtaining examples
of Lie $\mathfrak{g}$-foliations on compact manifolds when $\mathfrak{g}$
is compact:

\begin{Example} {\rm  Let $G$ be a connected compact Lie group, $N$ a compact
manifold and $h: \pi_1(N) \to G$ a homomorphism of the fundamental group of $N$ in $G$.
 Denote by $\tilde{N}$ the universal covering of $N$
and consider the diagonal action of $\pi_1(N)$ on $\tilde{N} \times G$ defined by
 $\alpha(x,g) \mapsto ((x)\alpha, h(\alpha)^{-1} g)$. Such action is properly discontinuous
 and this allow us to give a smooth structure which turns $E =  \left(\tilde{N} \times
 G\right)/ \pi_1(N)$ a (compact) manifold.  Moreover, we can define
 a foliation $\mathcal{F}$ on $E$ which is transverse to the fibers
 of $E$. Such foliation admits a Lie transverse structure of model $G$ and thus is a Lie $\mathfrak{g}$-foliation.}
\end{Example}

The result below is contained in Fedida's work:

{\bf \flushleft  Fedida's Theorem (\cite{Fedida}).}\label{fedida}{
\it Let $M$ be a compact manifold and $G$ a connected Lie group with
Lie algebra $\mathfrak{g}$. If $\mathcal{F}=\mathcal{F}_{\Omega}$ is
a Lie $\mathfrak{g}$-foliation on $M$ with global holonomy $\Gamma =
\Gamma(\Omega,G)$, then there is a locally trivial fibration $P: M
\to \overline{\Gamma}\setminus G$ such that each fiber
$\mathcal{F}$-saturated subset of $M$. In particular,
$\overline{\Gamma} \setminus G $ is compact.}

\begin{Remark}\label{simplesmente}{\rm When $G$ is a simply-connected Lie group, there exists a more
complete version of the Fedida's Theorem. In this case, the
fibration $P: M \to \overline{\Gamma}\setminus G$ has connected
fibers and $\mathcal{F}$ induces on each fiber of $P$ a minimal Lie
$\mathfrak{h}$-foliation, where $\mathfrak{h}$ is the Lie algebra of
$\overline{\Gamma}$.}
\end{Remark}

\section{Foliations invariant under transverse actions}

The Lie $\mathfrak{g}$-foliation in Example
2  has a special property: there is a smooth action $\phi$ on
the ambient manifold which preserve the leaves of $\mathcal{F}$ and whose orbits are transversal to $\mathcal{F}$.
The same is true for the foliation $\mathcal{G}$ defined on $M
\times G$ in Darboux-Fedida's Theorem. They are part of a special class
of foliations called {\it $G$-i.u.t.a.} foliations (\cite{behague}), which we describe more formally: Let $G$
be a $q$-dimensional connected  Lie group and $\mathcal{G}$ a codimension $q$ foliation
on a smooth manifold $E$. We say that $\mathcal{G}$ is  a
{\it $G$-i.u.t.a.} foliation (or {\it
invariant under transverse action of $G$}) if there is a smooth locally free action $\phi$ of $G$ on $E$ whose orbits
are transverse to $\mathcal{G}$ and such that, for each $g \in
G$, $\phi^g$ take leaves of $\mathcal{G}$ diffeomorphically
over leaves of $\mathcal{G}$. In \cite{behague}, the authors construct a development for this class of foliations. In particular,
if $\mathfrak{g}$ is the Lie algebra of $G$, then a $\gt$
foliation is a Lie $\mathfrak{g}$-foliation. The converse is not necessarily true (see \cite{behague}).

\begin{Proposition}[Fedida's Theorem for $\gt$ foliations]\label{Fedidagt} Let $E$ be a smooth closed manifold let and
$\mathcal{G}$ be a $G$-i.u.t.a. foliation on $E$. Let $L$ be a leaf of
$\mathcal{G}$ and consider the subgroup $\Gamma = \{g \in G ; \
\phi^g(L) =L\}$. There is a locally trivial fibration with connected fibers $P:E \to G / \overline{\Gamma}$
such that each fiber is a $\mathcal{G}$-saturated subset of $E$. Moreover, $\mathcal{G}$
induces on the fiber containing $L$ a minimal Lie $\mathfrak{h}$-foliation,
where $\mathfrak{h}$ is the Lie algebra of de $\overline{\Gamma}$.
\end{Proposition}

\begin{proof}[Proof] Let $\phi$ be the $G$-action on $E$ which preserves $\mathcal{G}$. Let $p_1: G \times L \to G$
be the canonical projection. It follows from Theorem 1 of \cite{behague}
that the map
$$\phi|_{G \times L} : G \times L \to E$$
is a Galoisian covering map and that there is a unique $P: E \to G /
\overline{\Gamma}$ such that the diagram
\begin{equation} \label{diagrama1}\xymatrix{ {G \times L} \rc{\phi|_{G \times L}} \ar[d]_{p_1} &
 E \ar[d]^{P} \\  {G} \ar[r]^{\pi} & G/
\overline{\Gamma}}
\end{equation}
is commutative. We conclude that $P$ is a submersion and, since $E$ is compact,
$P$ is a locally trivial fibration. We affirm that $\overline{\Gamma} = H =  \{g \in
G; \ \phi^g(\overline{L})= \overline{L} \}$. In fact,
the closed set $P^{-1}([e])$
contains $\overline{L}$. Therefore, by Diagram \ref{diagrama1}, we have that $H \subset \overline{\Gamma}$.
Since $E$ is compact, if $g\in \overline{\Gamma}$, then $\phi^g(\overline{L}) \subset
\overline{L}$. On the other hand, if $g \in \overline{\Gamma}$, then $g^{-1} \in \overline{\Gamma}$. Thus
$\phi^{g^{-1}}(\overline{L}) \subset \overline{L}$ and $\overline{\Gamma} \subset H$. Let $L_1$ be a leaf of $\mathcal{G}$
contained in $P^{-1}([e])$. There is $g \in G$ such that $\phi^g(L)= L_1$.
Since $L_1 \subset P^{-1}([e])$, we have that $g \in \overline{\Gamma}$.
Therefore, $L_1 \subset \phi^g(\overline{L}) =\overline{L}$,
and we conclude that $P^{-1}([e])= \overline{L}$. Since $\phi$
preserves $\mathcal{G}$, we have $\overline{L} =
\phi^g(\overline{L})= \overline{\phi^g(L)} = \overline{L_1}$.
Therefore, $\mathcal{G}$ induces on $\overline{L}$ a minimal foliation. If $H_0$ is the connected component of $H$
containing the identity, $\mathcal{G}$ induces on $\overline{L}$ a $H_0\textnormal{-i.u.t.a.}$ foliation and thus a
Lie $\mathfrak{h}$-foliation.
\end{proof}

\section{Lie foliations with compact Lie algebra}

In this section we study the structure of Lie $\mathfrak{g}$-foliations when $\mathfrak{g}$ is compact proving  Theorem~A.

\begin{proof}[Proof of the Theorem A] We can suppose that $M$ is orientable,
taking the orientable double covering if necessary. Let
$\mathcal{G}$ be the foliation on $E = M \times G$ constructed in
Darboux-Fedida Theorem. If $p_1: M \times G \to M$ and $p_2: M
\times G \to G$ are the canonical projections and $M'$ is a leaf of
$\mathcal{G}$, we have that $p_1|_{M'}: M' \to M$ is a Galoisian
covering. The automorphism group of this covering is associated with
a homomorphism $h: \pi_1(M) \to G$ whose image is
$$\Gamma = \{g \in G; \ \phi^g(M') = M' \},$$
the global holonomy of $\mathcal{F}$. Moreover, the natural action $\phi$ on $E= M \times G$ by left tranlations
turns $\mathcal{G}$ a $\gt$ foliation. From Fedida's Theorem (for $\gt$ foliations), there is a fibration $P: E \to G /
\overline{\Gamma}$ such that the fiber $N = P^{-1}([e])$ is the closure of the leaf $M'$.
Moreover, we have that
\begin{enumerate} \item  $\overline{\Gamma} = H = \{g \in G; \ \phi^g(N)= N\}$;
                  \item the foliation on $N$  induced for $\mathcal{G}$ is a minimal $\mathfrak{h}$-Lie foliation,
                      where  $\mathfrak{h}$ is a Lie algebra of $\overline{\Gamma}$.
\end{enumerate}
By hypothesis, $\overline{\Gamma}_0$ is Abelian and, thus,
$\mathcal{G}$ induces on $N$ a minimal Lie $\R^r$-foliation
$\mathcal{T}$, for some $r \in \N$. now we divide the proof in two
cases: {\bf \flushleft  First case.} $H$ is connected.
\newline \newline
For our proposes, we could consider only the restriction of $\mathcal{G}$ on $N$. However, we will construct a perturbation of the
foliation $\mathcal{G}$ in $E$. For this, we adapt a beautiful argument
presented by E. Ghys in the proof of Theorem A of \cite{ghys}.
The homogeneous space $G / H$ is orientable (see
Prop. 5.15 of \cite{bro}). Since $E$ is orientable, each fiber of $P$
is orientable. Fix an orientation in $N$. Dado $[g] \in G/ H$,
 we can consider an orientation in $N_g = \phi^g(N)$ induced by $\phi^g$.
 Since $H$ is connected, such a orientation is consistent, that is,
 it does not depend on the choice of $g \in [g]$.
 Consider $E$ endowed with a Riemannian metric invariant by $\phi$ and on each fiber of $P$, consider
 the metric induced by $E$. Let $\omega_1, \dots, \omega_r$ be smooth linearly independent closed forms on $N$
 which define $\mathcal{T}$. Given $[g] \in G/H$, consider the fiber $N_{g} = \phi^g
(N)$ of $P$. We can use the diffeomorphism $\phi^g$ and the forms
$\omega_1, \dots, \omega_r$ to obtain smooth closed forms $\omega_{(1,g)}, \dots, \omega_{(r,g)}$
on $N_g$.  These forms are well defined, that is, they does not depend on the choice
of $g \in [g]$. In fact, if $l \in H$ and $f=
\phi^l|_N : N \to N$, then $f$ is isotopic to identity, because $H$ is connected conexo.
Therefore, for each $i \in \{1, \dots,
r\}$, there is a function $h_i : N \to \R$ such that $f_* \omega_i  -
\omega_i = dh_i$. On the other hand, since $f$ preserves $\mathcal{T}$,
we have that $h_i$ is a  {\it basic function} for $\mathcal{T}$ (see
\cite{mol}, pag 34). This implies that $h_i$ is constant along the leaves of
$\mathcal{T}$. Since $\mathcal{T}$ has dense leaves, $h_i$ is constant and thus $f_*(\omega_i)= \omega_i$
for each $i \in \{1, \dots, r\}$, and the assertion is proved. Since
$\phi$ preserves the leaves of $\mathcal{G}$, we have that
$\omega_{(1,g)}, \dots, \omega_{(r,g)}$ define the foliation
induced by $\mathcal{G}$ on $N_g$. From Tischler's Theorem, there are arbitrarily small
harmonic closed forms of degree one $u_1, \dots, u_r$ on $N$
such that
$$\alpha_1 =\left(\omega_1 + u_1\right) , \dots, \alpha_r = \left(\omega_r + u_r\right)$$
are linearly independent and define a $\R^r$-Lie foliation on $N$ having all leaves compact.
Using the diffeomorphism $\phi^g$ and the forms $u_1, \dots,
u_r$ we obtain De Rham cohomology classes
$$[u_{(1,g)}], \dots, [u_{(r,g)}]$$
on $N_g$. Since $H$ is connected, this can be done consistently.
For each $k > k_0$ and each $i \in \{1, \dots, r\}$, consider
the only harmonic form $u_{(i,g)}$ in the cohomology class $[u_{(i,g)}]$.
For each $g \in G$, we have that
\begin{equation} \label{pull}\phi^g_* \left(\omega_i + u_i\right) = \phi^g_* \left(\omega_i + u_{(i,e)}\right) =  \left(\omega_{(i,g)} +
u_{(i,g)}\right).\end{equation}
In fact, $\phi^g_*
\left(\omega_i\right) = \omega_{(i,g)}$ and $\phi^g$ takes the cohomology class of $u_{(i,e)}$ in the
cohomology class of $u_{(i,g)}$.
But there is a only harmonic representative in the class $[u_{(i,g)}]$ and $\phi^g$ is
a positive isometry. Therefore, $\phi^g_* \left( u_{(i,e)}\right) =
u_{(i,g)}$ and assertion is proved. The forms
$$\omega_{(1,g)} := \left(\omega_{(1,g)} + u_{(1,g)}\right) , \dots, \omega_{(r,g)}:=\left(\omega_{(r,g)} +
u_{(r,g)}\right)$$ define on $N_g$ a Lie $\R^r$-foliation having all leaves compact.
Since $\phi$ acts transitively on $G/H$, we have a new foliation $\tilde{\mathcal{G}}$ on $E$ satisfying:
\begin{itemize} \item  In the fiber $N_g$, the foliation $\tilde{\mathcal{G}}$ é defined by the
closed forms $\left(\omega_{(i,g)} + u_{(i,g)}\right)$,  para $i = 1, \dots, r$;
                \item $\tilde{\mathcal{G}}$ is invariant by $\phi$ and has all leaves compact.
\end{itemize}
Since $E$ is compact, if $u_1, \dots, u_r$ are sufficiently small, the foliation $\tilde{\mathcal{G}}$
remains transversal to $\phi$. If $M''$ is a leaf of $\tilde{\mathcal{G}}$, we have that
$$p_1|_{M''} : M'' \to M$$
is a finite Galoisian covering and the group
$$\tilde{\Gamma}= \{g \in G; \ \phi^g(M'')=M''\}$$
is finite. Moreover, since $N$ is compact, if $u_1, \dots, u_r$ are
sufficiently small, we have that
$$p_2|_{M''} : M'' \to G$$
is a submersion and, since $M''$ is compact, $p_2|_{M''}$ is a locally trivial fibration
. The diagram
\begin{equation}\label{diag2}\xymatrix{ {M''} \ar[r]^{p_1|_{M''}} \ar[d]_{p_2|_{M''}}  & {M}  \\
{G} & {}}\end{equation} provides a development to a Lie $\mathfrak{g}$-foliation
$\tilde{\mathcal{F}}$ on $M$ having all leaves compact. Taking the forms $u_1, \dots, u_r$ sufficiently small, the foliation $\tilde{\mathcal{F}}$ is $C^0$-close to $\mathcal{F}$.

{\bf \flushleft  Second case. } $H$ is not connected. \newline

In this case, we construct the following structures:

\begin{enumerate} \item Two finite coverings maps $\xi: \tilde{E} \to E$, $\eta: \tilde{M} \to M$ and
a locally trivial fibration $\tilde{p}_1 : \tilde{E} \to \tilde{M}$ such that the diagram
$$\xymatrix{ {\tilde{E}} \ar[r]^{\tilde{p}_1} \ar[d]_{\xi} & {\tilde{M}} \ar[d]^{\eta} \\ {E} \ar[r]^{p_1} &
{M}}$$ is commutative;
                  \item A free action $\tilde{\phi}$ of $G$ on $\tilde{E}$ such that
                  \begin{itemize}   \item $\xi \circ \tilde{\phi}^g \equiv \phi^g \circ \xi \textnormal{ for all
                  } g \in G$;
                                    \item $\tilde{p}_1(x) = \tilde{p}_1(y)$ if and only if $x$ and $y$
                                        are in the same orbit of $\tilde{\phi}$;
                                    \item $\tilde{\phi}$ turns $\tilde{p}_1$ a $G$-principal bundle;
                  \end{itemize}
                  \item A locally trivial fibration $\tilde{P} : \tilde{E} \to G / H_0$ with connected fibers such that
                  \begin{itemize} \item the restriction of $\xi$ to a fiber of $\tilde{P}$ defines a diffeomorphism with a
                  fiber of $P$;
                                  \item $\tilde{P}$ possesses a fiber $\tilde{N}$ which projects by $\xi$
                                      over $N$, for which
                                  $$ \{g \in G; \ \tilde{\phi}^g(\tilde{N}) = \tilde{N} \} = H_0 =\overline{\Gamma}_0.$$
                  \end{itemize}
\end{enumerate}

After obtaining such structures, we can consider the foliations $\xi^*(\mathcal{G})$
and $\eta^*(\mathcal{F})$ and proceed as in the previous case.
We go on to describe such structures. Fix the leaf
$N = P^{-1}([e])$ of $\mathcal{G}$. the action
\begin{eqnarray*} H_0 \times (G \times N ) &\to& (G \times N)
\\ \left(h_0,(g,p)) \right) &\mapsto&
 \left( g h_0^{-1}, \phi^{h_0}(p)\right)
\end{eqnarray*}
is proper and free, thus, there is a unique smooth structure of $H_0
\setminus (G \times N)$ which turs the canonical projection
$$\beta: (G \times N) \to H_0\setminus (G \times N)$$
a principal bundle. Consider the manifold $\tilde{E}: = H_0\setminus
(G \times N)$ endowed with this structure. The action
\begin{eqnarray*} \varphi: \big(H_0 \setminus H\big) \times \tilde{E}  &\to& \tilde{E} \\
\left((h), [(g,p)]\right) & \mapsto& \left[
\left(gh^{-1},\phi^h(p)\right) \right]
\end{eqnarray*}
is well-defined. In fact, if $\tilde{h} = h_0 h$, $\tilde{g} = g
h_0^{-1}$ and $\tilde{p} = \phi^{h_0}(p)$, then
\begin{eqnarray*}  \tilde{g} \tilde{h}^{-1} =  g h^{-1}
\underbrace{(h h_0^{-1} h^{-1})}_{l_0^{-1} \in H_0} h_0^{-1} =
gh^{-1}(h_0 l_0)^{-1}
\end{eqnarray*}
and
$$\phi^{\tilde{h}}(\tilde{p}) = \phi^{h_0 h h_0} (p)= \phi^{h_0 h h_0 h^{-1} h} (p) = \phi^{h_0 l_0} \circ
\phi^h(p).$$ Therefore, $\varphi$ is well-defined. Moreover,
$\varphi$ is free. Since $H_0 \setminus H$ is a finite group, there
is a unique smooth structure for $\varphi \setminus \tilde{E}$ which
turns
$$q: \tilde{E} \to \varphi \setminus \tilde{E}$$
a finite Galoisian map. The map
\begin{eqnarray*} \xi :\tilde{E} &\to& E \\ {[}(g,p){]} &\mapsto& \phi^{g}(p)
\end{eqnarray*}
is well-defined and a submersion, because $\phi = \xi \circ \beta: G
\times N \to E$. There is a unique map $f: \varphi \setminus
\tilde{E} \to E$ which which make the diagram bellow to commute.
$$\xymatrix{ {} & {\tilde{E}} \ar[ld]_q \ar[rd]^{\xi} & {} \\
 {\varphi \setminus \tilde{E}} \ar[rr]^{f} &{}& {E}}$$
and, thus, $f$ is smooth. Moreover, $f$ is bijective and this implies
that $\xi: \tilde{E} \to E$ is a finite covering map. Consider the free action
\begin{eqnarray*} H_0 \times N &\to& N \\ (h_0,p) &\mapsto& \phi^{h_0}(p)
\end{eqnarray*}
and the unique smooth structure of $H_0 \setminus N$ which turns the
canonical projection
$$i: N \to H_0 \setminus N$$
 a principal bundle.
Write $\tilde{M}:= H_0 \setminus N$. The action
\begin{eqnarray*} \psi: \left(H_0 \setminus H \right)
\times \tilde{M} &\to& \tilde{M} \\ ((h),[p]) &\mapsto&
\left[ \phi^h(p)\right]
\end{eqnarray*}
is well-defined and free. Consider the unique structure of $\psi
\setminus \tilde{M}$ which turns the projection
$$j:\tilde{M} \to \psi \setminus \tilde{M}$$
a fibration. Such a projection is a finite covering map, because
$H_0 \setminus H$ is finite. It follows from definition of $\phi$
that the map
\begin{eqnarray*} \eta: H_0 \setminus N = \tilde{M} &\to& M
\\ {[}p=(x,g){]} &\mapsto& x
\end{eqnarray*}
is well defined. Moreover, $\eta$ is smooth, because
$$p_1|_N = \eta \circ i: N \to M.$$
There is an unique smooth map $\theta: \psi \setminus \tilde{M} \to
M$ which make the diagram
$$\xymatrix{ {} & {\tilde{M}} \ar[ld]_j \ar[rd]^{\eta} & {} \\ {\psi \setminus \tilde{M}} \ar[rr]^{\theta} &{}& {M}}$$
to commute. This implies that $\theta$ is smooth. Moreover, $\theta$
is bijective and then $\eta$ is a finite covering map. The map
\begin{eqnarray*} \tilde{p}_1: \tilde{E} &\to& \tilde{M} \\ {[}(p,g){]} &\mapsto& [p]
\end{eqnarray*}
is a submersion, because the diagram
$$\xymatrix{ {G \times N} \ar[d] \rc{\beta} & {\tilde{E}} \ar[d]^{\tilde{p}_1}
\\{N} \ar[r]^i & {\tilde{M}}}$$
is commutative. Since $\tilde{E}$ is compact, we have that
$\tilde{p}_1$ is a locally trivial fibration which make the diagram of item 1 to commute and
we conclude item 1. By construction, the action
\begin{eqnarray*} \tilde{\phi}: G \times \tilde{E} &\to& \tilde{E}
\\ \left( l, \left[(g,p)\right]\right) &\mapsto& \left[(lg,p)\right]
\end{eqnarray*}
satisfies the conditions in item 2. The map
\begin{eqnarray*} \tilde{P}: \tilde{E} &\to& G/H_0 \\ {[}(g,p){]} &\mapsto& [g] \end{eqnarray*}
is well-defined and a submersion, because the diagram
$$\xymatrix{ {G \times N} \ar[d] \ar[r]^\beta & {\tilde{E}} \ar[d]^{\tilde{P}}
\\{G} \ar[r] & {G/H_0}}$$
is commutative. Since $\tilde{E}$ is compact, $\tilde{P}$ is a locally trivial fibration.
We have that
$$\tilde{N}= \left\{ {[}(h_0,p){]}; \  h_0 \in H_0  \textnormal{ and } p \in N \right\}= \tilde{P}^{-1}([e]) = \beta
(H_0 \times N)$$ is a connected set, because $N$ is connected. Moreover,
 $\xi|_{\tilde{N}} : \tilde{N} \to N$ is a diffeomorphism and
$$\{g \in G; \ \tilde{\phi}^g(\tilde{N})= \tilde{N} \} = H_0,$$
and this completes the proof.

\end{proof}

\begin{proof}[Proof of Corollary \ref{Corollary:i}]
In \cite{carcar}, Caron and Carrière prove that a minimal
one-dimensional Lie foliation on a closed manifold $M^n$ has Abelian
global holonomy (see \cite{car} for a detailed proof). Since the
foliation is minimal, the global holonomy is dense, and thus, such
foliation is a Lie $\R^{n-1}$-foliation. If $\mathcal{F}_{\Omega}$
is a (not necessarily minimal) one-dimensional Lie
$\mathfrak{g}$-foliation
 on $M$ and $G$ is the simply-connected Lie group having $\mathfrak{g}$ as Lie algebra, then, by Fedida Theorem,
$\overline{\Gamma (\Omega, G)}$ has Abelian Lie algebra. Since
$\mathfrak{g}$ é semisimple and compact, $G$ is compact and we
obtain Corollary \ref{Corollary:i} from Theorem \ref{Theorem:A}.
\end{proof}

Theorem 1 was obtained from a perturbation of the foliation $\gt$ on the compact manifold $E=M \times G$ constructed from original Lie foliation
on $M$. We can use the constructions developed during the proof of Theorem 1 to obtain a fact about $\gt$ foliations:

\begin{Proposition}\label{Propositiongt2} Let $E$ be a smooth orientable closed manifold and $\mathcal{G}$ a $G$-i.u.t.a. foliation on $E$,
where $G$ is a compact Lie group. Using the notations of Proposition
\ref{Fedidagt}, suppose that $\overline{\Gamma}_0$ is Abelian. Then,
up to a finite equivariant covering, $\mathcal{G}$ can be
$C^0$-approximated by a $G$-i.u.t.a. foliation having all leaves
compact.
\end{Proposition}

\begin{proof}[Proof] We can construct $\tilde{E}$ and $\tilde{\phi}$ as in the proof of Theorem \ref{Theorem:A}.
Taking a Riemannian metric on $\tilde{E}$ invariant by the action of $G$, we can proceed as in the first case in the proof of
Theorem \ref{Theorem:A} and to obtain a $\gt$ foliation on $\tilde{E}$ having all leaves compact.
\end{proof}

\section{Lie foliations on manifolds with amenable fundamental group}

Let $A$ be a (discrete) group and let $\mathcal{B}(A)$ be the vector
space of bounded real functions on $A$ endowed with the sup norm. A
{\it mean} on $A$ is a left-invariant linear functional
$\eta:\mathcal{B}(A) \to \R$ satisfying $\eta(f) \geq 0$ and
$\eta(1) =1$. The class of amenable groups is closed under
operations of taking subgroups, extensions, quotient groups and
direct limits. Solvable groups are amenable. In contrast, a group
having a non-cyclic free subgroup is not amenable. For more details,
see \cite{green} and \cite{day}. By Theorem 4.1 of \cite{tah}, a
finitely generated group having subexponential growth is amenable.

\begin{Lemma}\label{fecho}Let $G$ be a Lie group and $\Gamma <G$ a subgroup of $G$.
If $Q$ is a subgroup of $\Gamma$ such that $[\Gamma:Q] < \infty$, then
$[\overline{\Gamma} : \overline{Q}] \leq [\Gamma:Q] < \infty$.
\end{Lemma}

\begin{proof}[Proof]Write $\Gamma/ Q = \{ (y_1), \dots, (y_l) \}$ and suppose,
by contradiction, that $[\overline{\Gamma}:
\overline{Q}] > [\Gamma:Q] =l$. In this case, there are distinct elements $[g_1], \dots, [g_l],  [g_{l+1}]$
in $\overline{\Gamma} /
\overline{Q}$. For each $k \in \{1, \dots, l+1 \}$, there is a sequence $x^k_n$ in $\Gamma$
which converges to $g_k$. Taking subsequences if necessary, we can suppose that $(x^k_n) = (y_{j_k})
\in \Gamma/Q$. Therefore, for each $k \in \{1, \dots, l+1 \}$,
there is a sequence $q^k_n$ in $Q$ such that $x^k_n = y_{j_k} *
q^k_n$ and, thus, $(y_{j_k}^{-1} * g_k) \in \overline{Q}$. There are
$k_1, k_2 \in \{1, \dots, l+1\}$ such that $k_1 \ne k_2$ and
$y_{j_{k_1}} = y_{j_{k_2}} =y$. Then
$$(y^{-1} * g_{k_1})^{-1} * (y^{-1} * g_{k_2}) = g_{k_1}^{-1} * g_{k_2} \in \overline{Q}.$$
This is a contradiction, because $[g_{k_1}] \ne [g_{k_2}]$.
\end{proof}

\begin{Lemma}\label{tits}Let $G$ be a connected Lie group. If a finitely generated subgroup
$\Gamma< G$ is amenable, then there is a solvable subgroup $Q < \Gamma$ such that $[\Gamma : Q]
< \infty$.
\end{Lemma}

\begin{proof}[Proof] Let $Ad: G \to GL(\mathfrak{g}, \R)$ be the adjoint homomorphism.
Since $\Gamma$ is finitely generated, either $Ad(\Gamma)$ has a non-cyclic free
subgroup or $Ad(\Gamma)$ is virtually solvable (Tits Alternative \cite{Tits}).
Since $Ad(\Gamma)$ is an amenable group, the first option does not occur. Therefore, there is a solvable
subgroup $E$ of $Ad(\Gamma)$ such that $[Ad(\Gamma): E] < \infty$. Writing $T=
Ad^{-1}(E)$, we have that
$$\frac{T}{\ker(Ad|_T)} \simeq E.$$
Since $\ker(Ad|_T) < \ker Ad$ is Abelian ($G$ is connected), we
conclude that $T$ is solvable. Therefore, $Q =\Gamma \cap T$ is a
solvable group. Since $[Ad(\Gamma):E] < \infty$, we have that
$[\Gamma : Q] < \infty$ and the proof is concluded.
\end{proof}

\begin{Proposition}\label{Proposition:I} Let $M$ be a smooth closed manifold
with amenable fundamental group and let
$\mathcal{F}=\mathcal{F}_{\Omega}$ be a Lie $\mathfrak{g}$-foliation
on $M$. Let $G$ be a connected Lie group having $\mathfrak{g}$ as
Lie algebra and let $\Gamma=\Gamma(\Omega,G)$ be the global holonomy
of $\mathcal{F}$. The connected component of $\overline{\Gamma}$
containing the identity is a solvable Lie subgroup of $G$. In
particular, if $G$ is compact, then the connected component of
$\overline{\Gamma}$ containing the identity is Abelian.
\end{Proposition}

\begin{proof}[Proof.] Since $\pi_1(M)$ is amenable, $\Gamma$ is amenable.
Denote by $\overline{\Gamma}_0$ the connected component of $\overline{\Gamma}$
containing the identity. Let $Q$ be the solvable subgroup of  $\Gamma$ given by Lemma $\ref{tits}$.
Taking the closure of the derived series of $Q$, we conclude that $\overline{Q}$ is solvable.
Therefore $\overline{\Gamma}_0 \cap \overline{Q}$ is a closed subgroup of $\overline{\Gamma}_0$
and a solvable group. By Lemma \ref{fecho}, $[\overline{\Gamma}, \overline{Q} ] < \infty$ and this implies
that $[\overline{\Gamma}_0,\overline{\Gamma}_0 \cap
\overline{Q}]<\infty$. Since $\overline{\Gamma}_0$ is connected,
we conclude that $\overline{\Gamma}_0 = \overline{\Gamma}_0 \cap
\overline{Q}$ is solvable. Then, $\overline{\Gamma}_0$ is a connected solvable Lie group. If
 $G$ is compact, then $\overline{\Gamma}_0$ is compact and, thus, Abelian.
\end{proof}

\begin{proof}[Proof of Theorem \ref{Theorem:B}] It follows from Theorem \ref{Theorem:A} and Proposition \ref{Proposition:I}.
\end{proof}

Using the notations of Proposition \ref{Fedidagt}, we obtain

\begin{Proposition}\label{Propositiongt3} Let $E$ be a smooth orientable closed manifold with amenable fundamental group and let $G$ a compact Lie group. If $\mathcal{G}$  is a $G$-i.u.t.a. foliation on $E$ then, up to a finite equivariant covering, $\mathcal{G}$ can be
$C^0$-approximated by a $G$-i.u.t.a. foliation having all leaves compact.
\end{Proposition}

\begin{proof}[Proof]It follows from Proposition \ref{Propositiongt2} and Proposition \ref{Proposition:I}.
\end{proof}

\section{Subexponential growth and Lie foliations}

In this section, we study Lie foliations on manifolds with fundamental growth of subexponential growth.

\begin{Lemma}\label{nilpotente} Let $G$ be a connected Lie group and $\Gamma<G$ a finitely generated
subgroup of subexponential growth. There is a nilpotent subgroup $T< \Gamma$ such that
$[\overline{\Gamma}: \overline{T}] < \infty$.
\end{Lemma}

\begin{proof}[Proof] Since $\Gamma$ has subexponential growth, $\Gamma$
is amenable. Let $Q< \Gamma$ be the solvable subgroup given by Lemma $\ref{tits}$.
Since $[\Gamma: Q] < \infty$, $Q$ is a finitely generated solvable group
of subexponential growth. By Milnor-Wolf Theorem (\cite{milnor}, \cite{wolf}) $Q$
is virtually nilpotent, that is, there is a nilpotent subgroup $T$ of $Q$
$[Q:T] < \infty$. We have $[\Gamma:T] < \infty$ and, by Lemma \ref{fecho}, $[\overline{\Gamma}:\overline{T}] <
\infty$.
\end{proof}

A Lie group (not necessarily connected) is {\it compactly generated} if $G$
is generated by some symmetric and compact neighborhood $W$ of the identity of $G$.
If $W$ is such a neighborhood, we write
$$W^n = \underbrace{W \cdot W \cdots W \cdot W}_{n-\textnormal{times}}.$$
Using the Haar measure of $G$, we can consider the {\it volume growth } of $G$ and
this does not depend on the choice of $W$. In particular, we say that
 $G$ has { \it polynomial volume growth} if there are constants $c >0$ and $d \in \N^*$
 such that $\mu(W^n) \leq c
n^d, \ \ \forall n \in \N$. For more details, see \cite{guiv}.

\begin{Proposition}\label{Proposition:II} Let $M$ be a
smooth closed manifold with fundamental group of subexponential growth. If
$M$ admits a Lie $\mathfrak{g}$-foliation and
$G$ is a connected Lie group having $\mathfrak{g}$ as Lie algebra, then
$G$ has polynomial volume growth.
\end{Proposition}

\begin{proof}[Proof] Let $\mathcal{F}=\mathcal{F}_{\Omega}$ be such foliation and let
$\Gamma= \Gamma(\Omega, \mathcal{F})$ be the global holonomy of
$\mathcal{F}$. Using Fedida Fibration Theorem, we conclude that
$G/\overline{\Gamma}$ is compact. By \cite{mac}, $\overline{\Gamma}$
is compactly generated. It follows from Theorem I.4 of \cite{guiv}
that $G$ and $\overline{\Gamma}$ has the same volume growth. Let $T$
be
 the nilpotent subgroup of $\Gamma$ given by Lemma \ref{nilpotente}. Since $[\overline{\Gamma}:\overline{T}] <
\infty$, $\overline{T}$ is compactly generated and has the same volume growth of
of $\overline{\Gamma}$ and, thus, of $G$.
Taking the closure of the lower central series of $T$, we obtain a central series for $\overline{T}$ and, thus
 $\overline{T}$ is a nilpotent group. By Theorem II.4 of
 \cite{guiv}, $\overline{T}$ (and thus $G$) has polynomial volume growth.
\end{proof}

Given a Lie algebra $\mathfrak{g}$, there is a solvable ideal of $\mathfrak{g}$
that contains all the solvable ideals of $\mathfrak{g}$. This ideal is called {\it radical of $\mathfrak{g}$}
and we will denote it by $\mathfrak{r}(\mathfrak{g})$. We have that $\mathfrak{g}$ is solvable if and only if
$\mathfrak{r}(\mathfrak{g})=\mathfrak{g}$;
$\mathfrak{g}$ is {\it semisimple} if and only if
$\mathfrak{r}(\mathfrak{g})=0$, that is, if and only if
$\mathfrak{g}$ has no solvable ideal. If $G$ is a connected Lie group having
$\mathfrak{g}$ as Lie algebra, then the {\it radical of $G$} is the
Lie subgroup of $G$ having $\mathfrak{r}(\mathfrak{g})$ as Lie algebra.

\begin{Proposition}\label{Proposition:III} Let $M$ be
a smooth closed manifold with fundamental group of subexponential growth.
If $M$ admits a Lie $\mathfrak{g}$-foliation, and $G$ is a connected Lie group having
$\mathfrak{g}$ as Lie algebra, then $G/R$ is compact, where $R$ is the radical of $G$.
\end{Proposition}

\begin{proof}[Proof] By Proposition \ref{Proposition:II}, $G$ has
polynomial volume growth. By Theorem 6.39 of \cite{pater}, the adjoint representation of
$\mathfrak{g}$ has only pure imaginary eigenvalues (see (6.35) of \cite{pater}). It follows from 6.32 of \cite{pater}
that $G/R$ is compact.
\end{proof}

In order to obtain Theorem \ref{Theorem:C} in its most general form, we use the Levi-Malcev decomposition
for Lie algebras \cite{var}:

{\flushleft \bf Levi-Malcev's decomposition.}{ \it Let
$\mathfrak{g}$ be a Lie algebra and $\mathfrak{r}(\mathfrak{g}) \subset \mathfrak{g}$ its radical.
There is a semisimple subalgebra $\mathfrak{s} \subset
\mathfrak{g}$ such that  $\mathfrak{g}$ is a direct sum {\rm(}as vector space{\rm)} of $\mathfrak{s}$ and
$\mathfrak{r}(\mathfrak{g})$.}

\vglue.1in

A Lie subalgebra $\mathfrak{s}<\mathfrak{g}$ given by Levi-Malvec decomposition
is called {it Levi factor of $\mathfrak{g}$}. Let
$G$ be a connected Lie group having $\mathfrak{g}$ as Lie algebra.
If $S <G$ is the Lie subgroup of $G$ that has
$\mathfrak{s}$ as Lie algebra, then $S$ is called of {\it Levi
factor of $G$}.

\begin{Lemma}\label{Propositionmalcev} Let $\mathcal{F}$ be a Lie
$\mathfrak{g}$-foliation  on a smooth manifold $M$ and
$\mathfrak{s} < \mathfrak{g}$ a Levi factor of $\mathfrak{g}$.
There is a Lie $\mathfrak{s}$-foliation $\mathcal{S}$ on $M$.
Each leaf of $\mathcal{S}$ is a $\mathcal{F}$-saturated subset
of $M$.
\end{Lemma}

\begin{proof}[Proof] Consider a Levi-Malcev decomposition
 $\mathfrak{g}=\mathfrak{s} + \mathfrak{r}(\mathfrak{g})$.
Take a basis $\{X_1, \dots, X_s\}$ of $\mathfrak{s}$ and a basis
$\{X_{s+1},\dots,X_q\}$ of $\mathfrak{r}(\mathfrak{g})$. Let
$c_{ij}^k$ be the structure constants of $\beta = \{X_1, \dots,
X_s,X_{s+1},\dots,X_q\}$. Since $\mathfrak{r}(\mathfrak{g})$ is an ideal of $\mathfrak{g}$,
we have that
\begin{equation}\label{equation:decomposicao}\textnormal{if }1 \leq k \leq s \textnormal{ and either \emph{i} or \emph{j} are
are larger that \emph{s}, then }c_{ij}^k =0.
\end{equation}
Since $\mathcal{F}$ is a Lie $\mathfrak{g}$-foliation, there are
linearly independent forms of degree one
$$\omega_1, \dots,\omega_s,\omega_{s+1}, \dots, \omega_q$$
on $M$ which are tangent to $\mathcal{F}$
and satisfy the relations involving $\{c_{ij}^k\}$. By (\ref{equation:decomposicao}),
$\omega_1, \dots, \omega_s$ define a Lie $\mathfrak{s}$-foliation $\mathcal{S}$ on $M$.
\end{proof}

\begin{proof}[Proof of Theorem \ref{Theorem:C}] The first item is the Proposition \ref{Proposition:II}. Let $\mathfrak{s}$ be
a Levi factor of $\mathfrak{g}$. By Lemma \ref{Propositionmalcev}, there is a Lie $\mathfrak{s}$-foliation on $M$. If $S$
is a connected Lie group having $\mathfrak{s}$ as Lie algebra, then $S$ is semisimple. By Proposition \ref{Proposition:III}, $S$ is compact.
By Theorem \ref{Theorem:A}, there is a finite covering of $M$ which fibers over $S$ and, thus, we have item 2.
\end{proof}

\section{Extension to transversely homogeneous foliations}

\begin{proof}[Proof of Theorem \ref{Theorem:D}] Let $\mathfrak{g}$ and $\mathfrak{k}$ be the Lie algebras of $G$ and $K$, respectively.
By Theorem 2 and Corollary 3.2 of
\cite{blu}, there are $n$ smooth forms of degree one $\omega_1, \dots,
\omega_q, \omega_{q+1}, \dots, \omega_n$ on $M$ such that
\begin{itemize} \item[i.] $\omega_1, \dots, \omega_q$ are linearly independent at each point of $M$;
                \item[ii.] $d \omega_k = - \sum_{i <j} c^k_{ij} \omega_i \wedge
\omega_j$, where $\{c^k_{ij}\}$ are the structure constants of a basis $\beta= \{X_1, \dots, X_n\}$ of $\mathfrak{g}$ such that
$X_{q+1}, \dots, X_n$ is a basis of $\mathfrak{k}$.
\end{itemize}
Consider the form $\Omega = \sum_{l=1}^n \omega_l X_l$ of degree one on $M$ with values in $\mathfrak{g}$ and
 let $\alpha$ be the Maurer-Cartan form of $G$. Consider the form $\Theta$ on $M \times G$ with values in
$\mathfrak{g}$ given by
$$\Theta_{(x,g)}(u,v) = \Omega(x)u - \alpha(g)v,$$
as in Darboux-Fedida Theorem. The form $\Omega$ is singular in $M$, however,
the form $\Theta$ has rank $n$ in each point of $M
\times G$. Therefore, as in Darboux-Fedida's Theorem,
$\Theta$ induces a foliation $\gt$ $\mathcal{G}$
on $M \times G$ whose dimension is the dimension of $M$. Let $p_2:
M \times G \to G$ and $\pi: G \to G/K$ be the canonical projections. By
(i), if $M'$ is a leaf of $\mathcal{G}$, the map $\pi \circ
p_2|_{M'}: M' \to G/K$ is a submersion. Then, we can follow the proof of
Theorem \ref{Theorem:A} and use the Theorem \ref{Theorem:B} to obtain Theorem \ref{Theorem:D}.

\end{proof}

\begin{proof}[Proof of Corollary \ref{Corollary:ii}]By Corollary 2.1 of \cite{blu1}, the fundamental group of $M$ has polynomial
volume growth. Then, the first statement follows from Theorem \ref{Theorem:D}. In the case $q=1$, suppose that the global
 holonomy $\Gamma$ of $\mathcal{F}$ is not finite. Since $\overline{\Gamma}_0$ is Abelian, then $\overline{\Gamma}_0$ contains
  (a subgroup of $SO(3)$ isomorphic to) $SO(2)$. We have that
  $$SO(3) \supsetneqq \overline{\Gamma} \supset \overline{\Gamma}_0 \supset SO(2).$$
  Therefore, either $\overline{\Gamma} = SO(2)$ or $\overline{\Gamma} = O(2)$. It follows from construction in Theorem \ref{Theorem:A} that, in the first case, the finite covering mentioned
  is the identity and, in the second case, it has two folds.
\end{proof}

\section{A remark about Riemannian foliations}
This section is concerned with the study of Riemannian foliations
and contains results already found in \cite{PAMS}.  Let $M$ be a
smooth manifold and $\mathcal{F}$ a smooth foliation on $M$ of
codimension $q$. Denote by $\mathcal{X}(\mathcal{F})$ the set of
 smooth vector fields on $M$ which are tangent to $\mathcal{F}$.
 We say that a vector field $X$ on $M$ is
{\it foliated with respect to $\mathcal{F}$} if

$$[X,Y] \in \mathcal{X}(\mathcal{F}), \ \ \textnormal{ for all } Y \in
\mathcal{X}(\mathcal{F}).$$ If $N\mathcal{F}$ is the normal bundle
of $\mathcal{F}$ and $X$ is a foliated vector field, we can project
$X$ over $N \mathcal{F}$, obtaining a vector field $\overline{X}$
tangent to $N \mathcal{F}$. We say that $\overline{X}$ is the {\it
transversal vector field associated with $X$}. We say that a
foliation $\mathcal{F}$ is {\it transversely parallelizable} if its
normal bundle admits a trivialization by $q$ linearly independent
transversal foliated vector fields. Every Lie foliation is a
transversely parallelizable foliation, however, the reciprocal is
not necessarily true.

We say that a Riemannian metric $g$ on $M$ is {\it bundle-like} with respect to a foliation $\mathcal{F}$
if, given an open $U \subset M$ and foliated vector fields $X,Y$
on $U$ orthogonal to the leaves of $\mathcal{F}$, the function
$g(X,Y)$ is constant along the leaves of $\mathcal{F}|_U$. For all details,
see \cite{mol}. We say that $\mathcal{F}$ is a
{\it Riemannian foliation} if there exists a bundle-like metric with respect to $\mathcal{F}$.
A deep result due to Molino describes the structure of transversely parallelizable and
 Riemannian foliations and proves that the study of such foliations is closely related
 to the study of Lie foliations:

\vglue.1in
\noindent{\bf Molino's structure theorem  (\cite{Go,mol})}. {\it  Let $M$ be a connected closed manifold and let $\mathcal{F}$ be a foliation on $M$.
\begin{enumerate}  \item If $\mathcal{F}$ is transversally parallelizable, then the closures of the leaves of $\mathcal{F}$ are the fibers of a locally trivial basic fibration $\pi: M \to B$. The restriction of $\mathcal{F}$ to each fiber defines a Lie $\mathfrak{g}$-foliation with dense leaves. The Lie algebra $\mathfrak{g}$ is an algebraic invariant of $\mathcal{F}$ called {\it structural Lie algebra of $\mathcal{F}$}.
                   \item If $\mathcal{F}$ is a transversally orientable Riemannian foliation and $\tau : M' \to M$ is the principal $SO(q)$-bundle of orthonormal frames transverse to $\mathcal{F}$, then $\mathcal{F}$ lifts to a transversally parallelizable foliation $\mathcal{F}'$ on $M'$ such that
                   \begin{itemize} \item $\dim{\mathcal{F}'}= \dim{\mathcal{F}}$;
                                   \item $ \mathcal{F}'$ is invariant by the action of $SO(q)$.
                   \end{itemize}
\end{enumerate}
The structural Lie algebra of $\mathcal{F}'$ is an invariant of $\mathcal{F}$ called structural Lie algebra of $\mathcal{F}$.} \\

Based on Molino Theorem, we obtain the following:

\begin{Proposition}\label{Proposition:IV} Let $M$ be a smooth closed manifold
and let $\mathcal{F}$ be a Riemannian foliation on $M$.
If the fundamental group of $M$ is amenable, then the structural Lie algebra of $\mathcal{F}$
is solvable.
\end{Proposition}

\begin{proof}[Proof] Let
$$\xymatrix{ {\cdots} \ar[r] &  {\pi_2(M)} \rc{\phi}
& {\pi_1(SO(q))} \rc{\varphi} & {\pi_1(M')} \rc{\psi} & {\pi_1(M)} \ar[r] & {\cdots} }$$
be the homotopy exact sequence associated with the principal
$SO(q)$-bundle of orthonormal frames transverse to $\mathcal{F}$. We
have that $\ker{\psi} = {\rm Im} \varphi$ is Abelian and thus
amenable. Moreover, ${\rm Im}\psi$ is a subgroup of $\pi_1(M)$ and
thus it is amenable. Since
$$\pi_1(M')/ \ker{\psi} \simeq {\rm Im} \psi,$$
we conclude that $\pi_1(M')$ is amenable, because amenable  groups
are closed to extensions. Let $\pi: M' \to B$ the locally trivial
basic fibration of $\mathcal{F}'$ and let $N= \pi^{-1}(y)$ be a
fibre of $\pi$. From homotopy exact sequence
$$\xymatrix{ {\cdots} \ar[r] &  {\pi_2(B)} \rc{\phi} & {\pi_1(N)} \rc{\varphi} & {\pi_1(M')} \rc{\psi} & {\pi_1(B)} \ar[r] & {\cdots} }$$
of $\pi$, we have that ${\rm Im} \varphi = \ker \psi$ is a subgroup
of $\pi_1(M')$, so it is amenable. Moreover, ${\rm Im} \phi =
\ker{\varphi}$ is an Abelian group, so it is amenable. Since
$$ \pi_1(N) / \ker \varphi \simeq {\rm Im} \varphi,$$
we conclude that $\pi_1(N)$ is amenable. From Molino Theorem,
$\mathcal{F}$ induces a minimal Lie $\mathfrak{g}$-foliation on $N$.
Let $G$ be a connected lie group having $\mathfrak{g}$ as Lie
algebra.  By Fedida theorem, the global holonomy of $\mathcal{F}$ is
dense in $G$. It follows from Proposition \ref{Proposition:I} that
the structural Lie algebra of $\mathcal{F}$ is solvable.

\end{proof}

\bibliographystyle{amsalpha}

\begin{tabular}{ll}
Marcelo Tavares:  marcelotrls@gmail.com\\
Instituto de  Matem\'atica -  Universidade Federal do Rio de Janeiro\\
Rio de Janeiro - RJ,   Caixa Postal 68530\\
21.945-970 Rio de Janeiro-RJ \\
BRAZIL
\end{tabular}

\end{document}